\documentclass[a4paper,leqno]{article}
\usepackage{amssymb,amsmath,amsfonts,amsthm,mathrsfs,graphicx,color}
\newtheorem{theorem}{Theorem}[section]

\newtheorem{lemma}[theorem]{Lemma}
\newtheorem{proposition}[theorem]{Proposition}

\newtheorem{assumption}[theorem]{Assumption}

\newcommand{\Vint}{\iint_{Q}}
\newcommand{\Bint}{\int_{t_0}^{T^{\prime}}\hspace{-6.5pt}  \int_\Omega}
\newcommand{\bint}{\int_{t_0}^T\hspace{-6pt}  \int_{\Gamma_0}}
\begin{document}
\title
      {Inverse Conductivity Problem for a Parabolic
 Equation using a Carleman Estimate with One Observation }





\maketitle

\centerline{\scshape Patricia Gaitan }
\medskip
{\footnotesize
 \centerline{Laboratoire d'Analyse, Topologie, Probabilit\'{e}s}
   \centerline{CNRS UMR 6632, Marseille, France and Universit\'e de la M\'editerrann\'ee}
} 

\bigskip


\begin{abstract}
For the heat equation in a bounded domain we give a
 stability result for a smooth diffusion coefficient.
The key ingredients are a global Carleman-type estimate, a Poincar\'e-type estimate
and an energy estimate with a single observation acting on a part of the boundary.
\end{abstract}

\section{Introduction}
This paper is devoted to the identification 
of the diffusion coefficient 
in the heat equation using the least number of 
observations as possible.\\
Let $\Omega \subset \mathbb{R}^n$ be a bounded domain of $\mathbb{R}^n$
with $n \le 3$, (the assumption $n \le 3$ is necessary in order to obtain the appropriate regularity
for the solution using classical Sobolev embedding, see Brezis \cite{B:83}).
We denote $\Gamma=\partial \Omega$
assumed to be
of class $\mathcal{C}^1$. 
We denote by $\nu$ the outward unit normal to $\Omega$ on $\Gamma=\partial \Omega$. 
Let $T>0$ and $t_0 \in (0,T)$. We shall use the following notations $Q_0=\Omega \times (0,T)$,
$Q=\Omega \times (t_0,T)$, $\Sigma=\Gamma \times (t_0,T)$  
and $\Sigma_0=\Gamma \times (0,T)$.
We consider the following heat equation:
\begin{equation}
\left \{ \begin{array}{lll}
 \label{systq}
  \partial_t q= \nabla \cdot (c(x) \nabla q) & \mbox{in} & Q_0,\\
   q(t,x)=g(t,x) & \mbox{on} & \Sigma_0 ,\\
   q(0,x)=q_0  & \mbox{in} & \Omega.
\end{array}\right.
\end{equation}
Our problem can be stated as follows:\\
{\bf Inverse Problem}\\
Is it possible to determine the coefficient $c(x)$ from the following measurements:
$$\partial_{\nu}(\partial_t q)_{|(t_0,T) \times \Gamma_0}\ \mbox{ and }
\nabla(\Delta q(T^{\prime}, \cdot)),\Delta q(T^{\prime}, \cdot), \ q(T^{\prime}, \cdot)
\ \mbox{ in } \Omega \mbox{ for }  T^{\prime}=\frac{t_0+T}{2},$$
where $\Gamma_0$ is a part of the boundary $\Gamma$ of $\Omega$ ?
\\ \\
\noindent
Let $q$ (resp. $\widetilde{q}$) be solution of (\ref{systq}) associated to
($c$, $g$, $q_0$) (resp. ($\widetilde{c}$,$g$, $q_0$)),
we assume
\begin{assumption}\label{reg-coeff}
\begin{itemize}
\item
$q_0$ belongs to $H^4(\Omega))$and $g$
is sufficiently regular (e.g. $\exists \;\epsilon>0$ such that $g \in H^1(0,T,H^{3/2+\varepsilon}(\partial \Omega))
\cap  H^2(0,T,H^{5/2+\varepsilon}(\partial \Omega))$)
\item
$ c,\ \widetilde{c} \in \mathcal{C}^3(\Omega)$,
\item
There exist a constant $r>0$, such that $q_0 \geq r$ and $g \geq r$.
\end{itemize}
\end{assumption} 
Note that the first item of the previous assumptions implies
that (\ref{systq}) admits a
solution in $H^1(t_0,T, H^{2}(\Omega))$ (see Lions \cite{L:68}).
We will later use this regularity result.
The two last items allows us to state that the function $u$
satisfies $| \Delta q(x,T')| \ge r>0$ and $| \nabla q(x,T')| \ge r>0$ in $\Omega$ 
(see Pazy \cite{P:83}, Benabdallah, Gaitan and Le Rousseau \cite{BGLR:07}).\\ \noindent
 We assume that we can measure both the normal flux $\partial_{\nu}(\partial_t q)$ on $\Gamma_0 \subset \partial \Omega$ 
in the time interval $(t_0, T)$ for some $t_0 \in (0, T)$ and
$\nabla(\Delta q)$, $\Delta q$ and $\nabla q$ at time $T' \in (t_0, T)$.\\
\noindent
Our main result is a stability result for the coefficient $c(x)$: \\ \noindent
For $q_0$ in $H^2(\Omega)$ there exists a constant
$C=C(\Omega, \Gamma, t_0, T, r) >0$
such that
\begin{eqnarray*}  
    |c-\widetilde{c}|^2_{H^1_0(\Omega)} & \leq &
    C |\partial_{\nu}(\partial_t q)- \partial_{\nu}(\partial_t \widetilde{q})|^2_{L^2((t_0,T) \times \Gamma_0)}\\
    &+& C |\nabla(\Delta q (T', \cdot))- \nabla(\Delta \widetilde{q}(T', \cdot))|^2_{L^2(\Omega)}\\ 
    & + & C |\Delta q (T', \cdot)- \Delta \widetilde{q}(T', \cdot)|^2_{L^2(\Omega)}   
    +C |\nabla q (T', \cdot)- \nabla \widetilde{q}(T', \cdot)|^2_{L^2(\Omega)}.
\end{eqnarray*}
The key ingredients to this stability result are a global Carleman-type estimate, 
a Poincar\'e-type estimate and an energy estimate.
We use the classical Carleman estimate with one observation on the boundary for
the heat equation obtained in Fernandez-Cara and Guerrero \cite{FG:05}, 
Fursikov and Imanuvilov \cite{FI:96}. Following the method developed by
Imanuvilov, Isakov and Yamamoto for the Lam\'e system in Imanuvilov, Isakov and Yamamoto \cite{IIY:03}, we give a
Poincar\'e-type estimate. Then, we prove an energy estimate. Such energy estimate
has been proved in Lasiecka, Triggiani ang Zhang \cite{LTZ:03} for the Schr\"{o}dinger operator in a
bounded domain in order to obtain a controllability result and in Cristofol, Cardoulis and Gaitan \cite{CCG:06}
for the Schr\"{o}dinger operator in a unbounded domain in order to obtain 
a stability result.
Then using these estimates, we give
a stability and uniqueness result for the diffusion coefficient $c(x)$.
In the perspective of numerical reconstruction,
such problems are ill-posed and stability results are thus of importance.\\
In the stationnary case, the inverse conductivity problem has been studied
by several authors. There are different approaches. For the two dimensional
case, Nachman \cite{N:95} proved an uniqueness result for the diffusion coefficient
$c \in C^2(\overline{\Omega})$ and Astala and P\"{a}iv\"{a}rinta \cite{AP:06}
for $c \in L^{\infty}(\Omega)$ with many measurements from the whole boundary.
In the three dimensional case, with the use of complex exponentially solutions,
Faddeev \cite{F:65}, Calderon \cite{C:80}, Sylvester and Uhlmann \cite{SU:87}
showed uniqueness for the diffusion coefficient.\\
There are few results on Lipschitz stability for parabolic equations, we can cite
Imanuvilov and Yamamoto \cite{IY:98}, Benabdallah, Gaitan and Le Rousseau \cite{BGLR:07}.
In \cite{BGLR:07}, the authors prove a Lipschitz stability result for the
determination of a piecewise-constant diffusion coefficient. For smooth coefficients
in the principal part of a parabolic equation, Yuan and Yamamoto \cite{YY:06} give a Lipschitz 
stability result with multiple observations. This paper is an improvement
of the simple case in \cite{YY:06} where we consider that the diffusion coefficient is
a real valued function and not a $n\times n$-matrix.
Indeed, in this case, with the method developped by \cite{YY:06},
they need two observations in order to obtain an estimation of the
$H^1$-norm of the diffusion coefficient. In this case, we need
only one observation.\\
Our paper is organized as follows. 
In Section 2, we recall the global Carleman estimate 
for (\ref{systq}) 
with one observation on the boundary.
Then we prove a Poincar\'e-type estimate for the coefficient $c(x)$
and an energy estimate. 
In Section 3, using the previous results, we establish
a stability estimate for the coefficient $c(x)$ 
when one of the solutions 
$\widetilde{q}$ is in a particular class of solutions with some 
regularity and "positivity" properties.
%
\section{Some Usefull Estimates}
 \subsection{Global Carleman Estimate}
\vskip 0.3cm

We recall here a Carleman-type estimate with a single observation 
acting on a part $\Gamma_0$ of the boundary $\Gamma$ of $\Omega$ 
in the right-hand side of the estimate (see \cite{FG:05}), \cite{FI:96}. Let us
introduce the following notations:
 
let $\widetilde{\beta}$ be a 
$\mathcal {C}^4(\overline{\Omega})$ positive function such that there exists a positive
constant $C_0$ which satisfies 
\begin{assumption}
\label{funct-beta}

$|\nabla \widetilde{\beta}| \geq C_0>0 \;\;\mbox{ in } \;\;\Omega, \;\;
{\partial}_{\nu} {\widetilde{\beta}}\leq 0\;\;\mbox{on}\;\;\Gamma \setminus \Gamma_0$,

\end{assumption}

Then, we define $\beta= \widetilde{\beta}+K$ with
$K= m \|\widetilde{\beta}\|_{\infty}$ and $m>1$. For $\lambda> 0$ and $t \in
(t_0,T)$, we define the weight functions
$$
  \label{wf}
  \varphi(x,t)=\frac{e^{\lambda \beta(x)}}{(t-t_0)(T-t)},
  \quad \quad \eta(x,t)=\frac{e^{2\lambda K} -e^{\lambda
  \beta(x)}}{(t-t_0)(T-t)}.
$$
If we set $\psi=e^{-s \eta}q$, we also introduce the following
operators
\begin{eqnarray*}
  M_1\psi  & =& \nabla \cdot(c \nabla \psi) + s^{2}\lambda^{2} c|\nabla\beta|^2\varphi^{2}\psi + s(\partial_t{\eta})\psi, \\
  M_2\psi  & =& \partial_t \psi-+2 s \lambda \varphi c \nabla \beta . \nabla \psi -2
s \lambda^{2} \varphi c | \nabla \beta |^2 \psi.
\end{eqnarray*}
Then the following result holds (see \cite{FG:05}, \cite{FI:96})\\
\begin{theorem}
\label{th-Carl-Fur} There exist $\lambda_0=\lambda_0(\Omega,
\Gamma_0)\geq 1$, $s_0=s_0(\lambda_0, T)>1$ and a positive
constant $C=C(\Omega, \Gamma_0, T)$ such that, for any $\lambda \ge
\lambda_0$ and any $s \ge s_0 $, the following inequality holds:
\begin{eqnarray}
\label{Carl-Fur1}
    \|M_1(e^{-s \eta}q)\|^2_{L^2(Q)} + \|M_2(e^{-s \eta}q)\|^2_{L^2(Q)}\\
  +s \lambda^2 \int \hspace{-6.5pt} \int_Q e^{-2s \eta} \varphi |\nabla q|^2 \ d x\ d t
  +s^{3} \lambda^{4}\int \hspace{-6.5pt} \int_Q e^{-2s \eta} \varphi^{3} |q|^2\ d x\ d t  \nonumber\\
  \leq C  \left[
   s \lambda\int_{t_0}^T \hspace{-6.5pt}
\int_{\Gamma_0} e^{-2s \eta} \varphi |\partial_{\nu} q|^2\ d x\ d t
  +\int \hspace{-6.5pt} \int_Q e^{-2s \eta}\ | \partial_t q - \nabla \cdot (c \nabla q)|^2\ d x\ d t\right],\nonumber
\end{eqnarray}
for all $q \in H^1(t_0,T, H^2(\overline{\Omega}))$ with $q=0$ on
$\Sigma$.
\end{theorem} 

\subsection{Poincar\'e-type estimate}
\vskip 0.3cm
We consider the solutions $q$ and $\widetilde{q}$
to the following systems
\begin{equation}
\left \{ \begin{array}{lll}
  \label{syst-q}
    \partial_t q= \nabla \cdot (c(x) \nabla q) & \mbox{in} & Q_0,\\
    q(t,x)=g(t,x)& \mbox{on} & \Sigma_0,\\
    q(0,x)=q_0  & \mbox{in} & \Omega, 
\end{array}\right.
\end{equation}
and
\begin{equation}
\left \{ \begin{array}{lll}
  \label{syst-qtilde}
   \partial_t \widetilde{q}= \nabla \cdot (\widetilde{c}(x) \nabla \widetilde{q}) & \mbox{in} & Q_0,\\
    \widetilde{q}(t,x)=g(t,x)& \mbox{on} & \Sigma_0,\\
    \widetilde{q}(0,x)=q_0 & \mbox{in} & \Omega.
\end{array}\right.
\end{equation}
We set $u=q-\widetilde{q}$, $y=\partial_t u$ and $\gamma=c-\widetilde{c}$.
Then $y$ is solution to the following problem
\begin{equation}
\left \{ \begin{array}{lll}
  \label{syst-y}
     \partial_t y= \nabla \cdot (c(x) \nabla y)+\nabla \cdot (\gamma(x) \nabla (\partial_t \widetilde{q})) & \mbox{in} & Q_0,\\
       y(t,x)=0  & \mbox{on} & \Sigma_0,\\
    y(0,x)=\nabla \cdot (\gamma(x) \nabla (q_0(x))), & \mbox{in} & \Omega.
\end{array}\right.
\end{equation}

Note that with (\ref{syst-q}) and (\ref{syst-qtilde}) we can determine $y(T',x)$ and we obtain
\begin{equation} \label{CIT'}
y(T',x)=\nabla \cdot (\gamma(x) \nabla (\widetilde{q}(T',x)))
+\nabla \cdot (c(x) \nabla (u(T',x))).
\end{equation}

We use a lemma proved  in \cite{IIY:03} for Lam\'e system in bounded domains:

\begin{lemma}\label{IIY}
We consider the first order partial differential operator
$$ P_0g:=\nabla q_0 \cdot \nabla g$$
where $q_0$ satisfies 

$|\nabla \beta \cdot \nabla q_0|\ \neq 0$.\\

Then there exists positive constant, $s_1>0$ and $C=
C(\lambda,T')$
such that for all  $s \geq s_1$
$$s^2 \lambda^2 \int_{\Omega}  \ e^{-2s \eta(T')}\varphi(T') |g|^2\ dx \ dy
\leq C  \int_{\Omega} e^{-2s \eta(T')}\varphi^{-1}(T') \ |P_0 g|^2\ dx \ dy$$
with $T' = \frac{t_0+T}{2}$, $\eta(T'):=\eta(x,T')$, $\varphi(T'):=\varphi(x,T')$ and 
for $g \in H^1_0(\Omega)$.
\end{lemma}

We assume

\begin{assumption} \label{un}

$|\nabla \beta \cdot \nabla  \widetilde{q}(T')|\ \neq 0$, 

\end{assumption}
\vskip 0.3cm

\begin{proposition}\label{Poinc}
Let $ \widetilde{q}$ be solution of (\ref{syst-qtilde}). We assume that Assumption \ref{un}
are satisfied. Then there exists a positive constant
$C=C(T',\lambda)$ such that for $s$ large enough
( $s\geq s_1$), the following estimate
hold true

$$
s^2 \lambda^2 \int_{\Omega}  e^{-2s \eta(T')} \varphi(T')( |\nabla \gamma|^2\ + | \gamma|^2 ) \ dx 
\leq C  \int_{\Omega}  e^{-2s \eta(T')}\varphi^{-1}(T') \left(|\nabla y(T')|^2+|y(T')|^2\right)\ dx   $$
$$+ C  \int_{\Omega}  e^{-2s \eta(T')} \left(
|\nabla(\Delta u(T^{\prime})|^2 + |\Delta u(T^{\prime})|^2 + 
\nabla u(T^{\prime}|^2  \right)\ dx$$

for $\gamma \in H^2_0(\Omega)$.
\end{proposition}
\noindent
\begin{proof}
We are dealing with the following first order partial differential operators given by the equation (\ref{CIT'})
$$
P_0 (\gamma):= \sum_{i=1}^n \partial_{x_i} \widetilde{q}(T') \partial_{x_i} \gamma 
=y(T') - \gamma \Delta \widetilde{q}(T') - \nabla (c \nabla u)(T').
$$
We apply the lemma \ref{IIY} for this operator and  we can write :

\begin{equation} \label{IIY1}
\begin{array}{lcl}
s^2 \lambda^2 \int_{\Omega}  e^{-2s \eta(T')}\varphi(T') | \gamma|^2 \ dx  
\leq  C  \int_{\Omega}  e^{-2s \eta(T')}\varphi^{-1}(T') \left( |y(T')|^2 + |\gamma |^2 \right)\ dx   \\ [5mm]
 + C  \int_{\Omega}  e^{-2s \eta(T')} \left(
|\Delta u(T^{\prime})|^2 + 
|\nabla u(T^{\prime})|^2  \right)\ dx
\end{array}
\end{equation}

In the other hand, we use the $x_j$-derivative of the previous equation (\ref{CIT'}). 
So, for each $j$ we deal with the following first order partial differential operator :

$$
P_0 (\partial_{x_j} \gamma)  = \partial_{x_ j} (T') - \partial_{x_j}\gamma \Delta \widetilde{q}(T') - \gamma \Delta ( \partial_{x_j} \widetilde{q})(T') - \partial_{x_j}(\nabla (c \nabla u))(T').
$$

Then under assumption (\ref{un}):

\begin{eqnarray*}
s^2 \lambda^2 \int_{\Omega}  e^{-2s \eta(T')}\varphi(T') |\partial_{x_j}  \gamma|^2 \ dx  
\leq C  \int_{\Omega}  e^{-2s \eta(T')} \varphi^{-1}(T') |\partial_{x_j}  y(T')|^2 \ dx \\ 
 + C  \int_{\Omega}  e^{-2s \eta(T')}\varphi^{-1}(T') \left( |\partial_{x_j} \gamma|^2+ |\gamma|^2 +| \nabla \gamma|^2   + |\partial_{x_j} F|^2 \right) \ dx 
\end{eqnarray*}

So, adding for all $j$, we can write 

\begin{eqnarray}\label{IIY2}
s^2 \lambda^2 \int_{\Omega}  e^{-2s \eta(T')}\varphi(T') |\nabla  \gamma|^2 \ dx  
\leq C  \int_{\Omega}  e^{-2s \eta(T')} \varphi^{-1}(T')|\nabla  y(T')|^2 \ dx \\ \nonumber  
 + C  \int_{\Omega}  e^{-2s \eta(T')}\varphi^{-1}(T') \left( |\nabla \gamma|^2+ |\gamma|^2 +| \nabla(\Delta u(T')|^2 + |\Delta u(T')|^2 \right) \ dx 
\end{eqnarray}

Taking into account (\ref{IIY1}) and (\ref{IIY2}) and for $s$ large enough, we can conclude.
\end{proof}

\subsection{Estimation of $\int_{\Omega} e^{-2s\eta(T')}|y(T')|^2\ dx$}
\vskip 0.3cm
Let $T'= \frac{1}{2}(T+t_0)$ the point for which
$\Phi(t)=\frac{1}{(t-t_0)(T-t)}$ has its minimum value. 

We set $\psi = e^{-s \eta} y$. 
With the operator

\begin{equation}
\label{eq: M2} 
M_2\psi   = \partial_t \psi-+2 s \lambda \varphi c \nabla \beta . \nabla \psi -2
s \lambda^{2} \varphi c | \nabla \beta |^2 \psi,
\end{equation}
we introduce, following \cite{BP:02}, 

\begin{eqnarray*}
  \mathcal{I} = \int_{t_0}^{T'} \int_{\Omega} M_2 \psi\;\psi\;dx dt
\end{eqnarray*}
We have the following estimates.
 \begin{lemma} 
  \label{lemma2}
  Let $\lambda \geq \lambda_1$, $s\geq s_1$ and let $a,\ b,\ c,\ d\ \in \L^{\infty}(\Omega)$. 
Furthermore,
we assume that $u_0$, $v_0$ in $H^2(\Omega)$ and the
assumption (\ref{reg-coeff}) is satisfied. Then there exists a constant $C=C(\Omega, \omega, T)$ such that
  \begin{equation}
  \label{eq:lemma2}
    \int_{\Omega} 
    e^{{-2s \eta}(T',x)}\
    |y(T',x)|^2 \; dx
    \leq C \left[ \lambda^{1/2}
   \int_{t_0}^{T} \int_{\Gamma_0} e^{-2s \eta} \varphi | \partial_{\nu} y |^2\ d x\ d t \right. 
    \end{equation}
    $$\left. + s^{-1/2} \lambda^{-1/2}\int_{t_0}^{T} \int_{\Omega} 
    e^{-2s \eta} \left(|\gamma|^2 +| \nabla \gamma |^2 \right) d x\ d t\right].$$
\end{lemma}
\noindent
\begin{proof} 
 If we compute $\mathcal{I}$, we obtain :
 $$\int_{\Omega} e^{{-2s \eta}(T',x)}\ |y(T',x)|^2 \; dx=-2 \mathcal{I}$$
 $$-4 s \lambda \int_{t_0}^{T'} \int_{\Omega} \varphi\ c \nabla \beta \cdot \nabla \psi\ \psi \ d x\ d t
 -4s \lambda^2 \int_{t_0}^{T'} \int_{\Omega} \varphi\ c |\nabla \beta|^2 |\psi|^2 \ d x\ d t.$$
 Then with the Carleman estimate (\ref{Carl-Fur1}), we can estimate all the terms in the
 right hand side of the previous equality and we have
 \begin{eqnarray*}
 \int_{\Omega} e^{{-2s \eta}(T',x)}\ |y(T',x)|^2 \; dx \leq
 C s^{-3/2} \lambda^{-2}\left( \|M_2 \psi\|^2
 +s^3 \lambda^4 \Vint e^{-2 s \eta} \varphi^3|y|^2 \ d x\ d t \right)\\
 + C s^{-1} \lambda^{-1/2} \left(s \lambda \Vint e^{-2 s \eta} \varphi\ |\nabla y|^2 \ d x\ d t
+s^3 \lambda^4 \Vint e^{-2 s \eta} \varphi^3|y|^2 \ d x\ d t \right)\\
+C s^{-2} \lambda^{-2} \left( s^3 \lambda^4 \Vint e^{-2 s \eta} \varphi^3|y|^2 \ d x\ d t \right).  
 \end{eqnarray*}
 Finally, we obtain
$$\int_{\Omega} e^{{-2s \eta}(T',x)}\ |y(T',x)|^2 \; dx \leq
C \lambda^{1/2} \bint e^{-2 s \eta} \varphi\ |\partial_{\nu} y|^2 \ d \sigma \ d t
+Cs^{-1} \lambda^{-1/2}\Vint e^{-2 s \eta}  |f|^2 \ d x\ d t,$$
where $f=\nabla \cdot(\gamma \nabla \partial_t \widetilde{q})$.
We assume that $\widetilde{q}$ is sufficiently smooth in order to have
$\nabla \partial_t \widetilde{q}$ and $\Delta \partial_t \widetilde{q}$ in
$L^2(O,T,L^{\infty}(\Omega))$.\\
Moreover taking into account that
$e^{-2 s \eta(t)} \leq e^{-2 s \eta(T')}$, the proof of Lemma \ref{lemma2} is complete.
\end{proof}
\subsection{Estimation of $\int_{\Omega} e^{-2s\eta(T')}\varphi^{-1}(T')|\nabla y(T')|^2\ dx$}
\vskip 0.3cm
We introduce
\begin{equation}
\label{energy}
E(t)=\int_{\Omega} c \ \varphi^{-1}(x,t)e^{-2s \eta(x,t)}|\nabla y(x,t)|^{2} \ dx.
\end{equation}
In this section, we give an estimation for the energy $E(t)$ at $T'$.\\
\begin{theorem}
\label{Energy-Estimate}
We assume that Assumptions \ref{reg-coeff} are checked, then there exist 
$\lambda_1=\lambda_1(\Omega, \omega)\geq 1$, $s_1=s_1(\lambda_1, T)>1$ and a positive constant $C=C(\Omega,
\Gamma_0, C_0, r, T)$ such that, for any $\lambda \geq \lambda_1$ and
any $s \geq s_1 $, the following inequality holds:
\begin{equation}
\label{EE}
E(T') \leq C \left[
    s \lambda  \int_{t_0}^T \int_{\Gamma_0} e^{-2s \eta} \varphi | \partial_{\nu} y |^2\ d x\ d t 
   + s \Vint e^{-2s \eta} (|\gamma|^2+|\nabla \gamma|^2) \ d x \ d t \right],
   \end{equation}
 \end{theorem} 
\noindent 
\begin{proof}
We note $f=\nabla \cdot (\gamma(x) \nabla \partial_{t} \widetilde{q})$.\\
We multiply the first equation of (\ref{syst-y}) by $e^{-2s \eta} \nabla \cdot (c \nabla y)\varphi^{-1}$ and
integrate over $(t_0,T)\times \Omega$, we have :
\begin{equation}
\label{EE:1}
\Bint \varphi^{-1}e^{-2s \eta}\nabla \cdot (c \nabla y)\partial_{t}y \ dx \ dt =
\Bint \varphi^{-1}e^{-2s\eta}|\nabla \cdot (c \nabla y)|^{2} \ dx \ dt
\end{equation}
$$
+\Bint e^{-2s \eta}\varphi^{-1} \nabla \cdot (c \nabla y)f \ dx \ dt. 
$$
we denote $\displaystyle{A:=\Bint e^{-2s \eta}\varphi^{-1}\nabla \cdot (c \nabla y)\partial_{t}y\ dx \ dt}$.\\
Integrating by parts $A$ with respect to the space variable, we obtain
\begin{align}
\label{EE:2}
A=\Bint c\ e^{-2s \eta}\varphi^{-1}\nabla y\partial_{t}(\nabla y) \ dx \ dt
+ 2s\lambda \Bint c \ e^{-2s \eta}\nabla y\partial_{t}y \nabla \beta \ dx \ dt
\end{align}
$$-\lambda \Bint c \ e^{-2s \eta}\varphi^{-1}\nabla y\partial_{t}y \nabla \beta \ dx \ dt.$$
Observe that
$$e^{-s \eta}\varphi^{-\frac{1}{2}}\partial_{t}(\nabla y)=\partial_{t}(e^{-s \eta}\varphi^{-\frac{1}{2}}\nabla y)+s e^{-s
\eta}\varphi^{-\frac{1}{2}}\partial_{t}\eta \nabla y + 
+ \displaystyle{\frac{1}{2}} e^{-s \eta}\partial_{t}\varphi \varphi^{-\frac{3}{2}}\nabla y.$$
Hence, the first integral of the right-hand side of (\ref{EE:2}) can be written as

$$\Bint c \ e^{-2s \eta}\varphi^{-1}\nabla y\partial_{t}(\nabla y) \ dx \ dt=\Bint c\ e^{-s \eta}\varphi^{-\frac{1}{2}}\nabla
y\partial_{t}(\nabla y)e^{-s \eta}\varphi^{-\frac{1}{2}} \ dx \ dt$$
$$=\Bint c\ e^{-s \eta}\varphi^{-\frac{1}{2}}\nabla y \partial_{t}(e^{-s \eta}\varphi^{-\frac{1}{2}}\nabla y) \ dx \ dt 
+s\Bint c\ e^{-2s \eta}\varphi^{-1}|\nabla y|^{2}\partial_{t}\eta \ dx \ dt$$
\begin{equation}
\label{EE:3}
+\displaystyle{\frac{1}{2}}\Bint c\ e^{-2s \eta}\varphi^{-2}|\nabla y|^{2}\partial_{t}\varphi \ dx \ dt.
\end{equation}
Using an integration by parts with respect the time variable, the first term of (\ref{EE:3}) is exactly equal to
$\displaystyle{\frac{1}{2}}E(T')$, since $E(t_0)=0$.
Therfore, the equations (\ref{EE:1}), (\ref{EE:2}) and (\ref{EE:3}) yield
$$
E(T')=
-2s\Bint c\ e^{-2s \eta}\varphi^{-1}|\nabla y|^{2}\partial_{t}\eta \ dx \ dt    
-\Bint c\ e^{-2s \eta}\varphi^{-2}|\nabla y|^{2}\partial_{t}\varphi \ dx \ dt$$
$$-4s\lambda \Bint c \ e^{-2s \eta}\nabla y\partial_{t}y \nabla \beta \ dx \ dt    
+2\lambda \Bint c \ e^{-2s \eta}\varphi^{-1}\nabla y\partial_{t}y \nabla \beta \ dx \ dt$$
$$+2\Bint \varphi^{-1}e^{-2s\eta}|\nabla \cdot (c \nabla y)|^{2} \ dx \ dt     
+2\Bint e^{-2s \eta}\varphi^{-1} \nabla \cdot (c \nabla y)f \ dx \ dt $$
\begin{equation}
\label{EE:4}
=I_1+I_2+I_3+I_4+I_5+I_6. 
\end{equation}
Now, in order to obtain an estimation to $E(T')$, we must estimate all the integrals $I_i, 1 \leq i \leq 6.$\\
Using the fact that $|\partial_t \eta| \leq C(\Omega, \omega) T \varphi^2$,
we obtain, in first step, for the integral $I_1$, the following estimation
\begin{eqnarray*}
|I_1| \leq  C s\Bint c\ e^{-2s \eta}\varphi|\nabla y|^{2} \ dx \ dt\\
\leq  C \lambda^{-2} \left[s \lambda^{2} \Vint e^{-2s \eta}\varphi|\nabla y|^{2} \ dx \ dt \right].
\end{eqnarray*} 
In a second step, the Carleman estimate yields
$$ 
|I_1| \leq  C \lambda^{-2} \left[ s \lambda  \bint e^{-2s \eta} \varphi
    | \partial_{\nu} y |^2\ d x\ d t
      + \Vint e^{-2s \eta} |f|^2 | \ d x \ d t \right],
$$
where $C$ is a generic constant which depends on $\Omega$, $\Gamma_0$, $c_{\max}$ and $T$.\\
As the same way, we have, for $I_2$, the following estimate 
$$
|I_2| \leq  C s^{-1}\lambda^{-2} \left[s \lambda  \bint e^{-2s \eta} \varphi
    | \partial_{\nu} y|^2\ d x\ d t
      + \Vint e^{-2s \eta} |f|^2 | \ d x \ d t \right].
$$
The last inequality holds throught the Carleman estimate and the following inequality
$$ |\partial_t \varphi| \leq C(\Omega, \Gamma_0) T^3 \displaystyle{\frac{\varphi^{3}}{4}}.$$
Using Young inequality, we estimate $I_3$.\\ 
We have
 \begin{eqnarray*}
  |I_3 | \leq  C s \left[s \lambda^{2} \Vint e^{-2s \eta} \varphi |\nabla y|^2 \ d x\ d t
+  s^{-1} \Vint e^{-2s \eta} \varphi^{-1} 
|\partial_t y|^2\ d x\ d t \right]\\
\leq  
C s \left[s \lambda  \bint e^{-2s \eta} \varphi
    | \partial_{\nu} y |^2\ d x\ d t
      + \Vint e^{-2s \eta} |f|^2 | \ d x \ d t \right],
\end{eqnarray*}
For the integral $I_4$, we have  
\begin{eqnarray*} 
|I_4| \leq  C  \left[s \lambda^{2} \Vint e^{-2s \eta} \varphi^{-1}  |\nabla y|^2 \ d x\ d t
+  s^{-1} \Vint e^{-2s \eta} \varphi^{-1} 
|\partial_t y|^2\ d x\ d t \right]\\
\leq  
C \left[s \lambda  \bint e^{-2s \eta} \varphi
    | \partial_{\nu} |^2\ d x\ d t
      + \Vint e^{-2s \eta} |f|^2 | \ d x \ d t \right],\\
\end{eqnarray*}
where we have used, for the term containing $|\nabla y|^2$, the following estimate
$$\varphi^{-1} \leq C(\Omega, \omega) T^4 \displaystyle{\frac{\varphi}{16}}.$$
we have immediatly the following estimate for $I_5$
\begin{eqnarray*}
|I_5| \leq C s \left[s \lambda  \bint e^{-2s \eta} \varphi
    | \partial_{\nu} y |^2\ d x\ d t
      + \Vint e^{-2s \eta} |f|^2 | \ d x \ d t \right].\\
\end{eqnarray*}

Finally, for the last term $I_6$, we have 
\begin{eqnarray*}
|I_6 | \leq  C \left[s^{-1} \Vint e^{-2s \eta} \varphi^{-2} |\nabla \cdot (c \nabla y)|^{2} \ d x\ d t
+  s \Vint e^{-2s \eta} |f|^2 \ d x\ d t \right]\\
\leq  
C \left[s \lambda  \bint e^{-2s \eta} \varphi
    | \partial_{\nu} y |^2\ d x\ d t
      + s \Vint e^{-2s \eta} |f|^2 | \ d x \ d t \right].\\
\end{eqnarray*}
The last inequality holds using the following estimate
$$\varphi^{-2} \leq C(\Omega, \omega) T^2 \displaystyle{\frac{\varphi^{-1}}{4}}$$
If we come back to (\ref{EE:4}), using the estimations of $I_i, 1 \leq i \leq 6$ and expanding the term $f$, 
this conlude the proof of Theorem \ref{Energy-Estimate}.
\end{proof}

\section{Stability Result}
\begin{theorem} \label{stab}
 Let $q$ and $\widetilde{q}$ be solutions of (\ref{syst-q}) and (\ref{syst-qtilde})
such that $c-\widetilde{c} \in H^2_0(\Omega)$.
We assume that  Assumptions \ref{reg-coeff} are satisfied.
  Then there exists a positive constant
  $C=C(\Omega, \Gamma_0,T)$ 
such that for $s$ and $\lambda$ large enough,
$$\int_{\Omega} \varphi(T') \ e^{-2s\eta(T')} (|c-\widetilde{c}|^2 +|\nabla (c-\widetilde{c}) |^2) \
dx \ dy 
\leq C  \int_{0}^T \int_{\Gamma_0} \varphi \ e^{-2s\eta} \partial_{\nu}
\beta \ |\partial_{\nu} (\partial_t q -\partial_t \widetilde{q})|^2 \ d\sigma \ dt$$
 $$+C  \int_{\Omega}  e^{-2s \eta(T')} \left(
|\nabla(\Delta u(T^{\prime})|^2 + |\Delta u(T^{\prime})|^2 + 
\nabla u(T^{\prime}|^2  \right)\ dx$$
 \end{theorem}
 \noindent
\begin{proof}
Using the estimates (\ref{EE}), (\ref{eq:lemma2}) and Proposition (\ref{Poinc}),
we obtain
\begin{eqnarray*}
s^2 \lambda^2 \int_{\Omega}  e^{-2s \eta(T')} \varphi(T')( |\nabla \gamma|^2\ + | \gamma|^2 ) \ dx 
\leq C  \int_{\Omega}  e^{-2s \eta(T')}\varphi^{-1}(T') \left(|\nabla y(T')|^2+|y(T')|^2\right)\ dx   \\
+ C  \int_{\Omega}  e^{-2s \eta(T')} \left(
|\nabla(\Delta u(T^{\prime})|^2 + |\Delta u(T^{\prime})|^2 + 
\nabla u(T^{\prime}|^2  \right)\ dx\\
\leq C \left[
    s \lambda  \int_{t_0}^T \int_{\Gamma_0} e^{-2s \eta} \varphi | \partial_{\nu} y |^2\ d x\ d t 
   + s \Vint e^{-2s \eta} (|\gamma|^2+|\nabla \gamma|^2) \ d x \ d t \right]\\
   + C \left[ \lambda^{1/2}
   \int_{t_0}^{T} \int_{\Gamma_0} e^{-2s \eta} \varphi | \partial_{\nu} y |^2\ d x\ d t 
     + s^{-1/2} \lambda^{-1/2}\int_{t_0}^{T} \int_{\Omega} 
    e^{-2s \eta} \left(|\gamma|^2 +| \nabla \gamma |^2 \right) d x\ d t\right]\\
    +C  \int_{\Omega}  e^{-2s \eta(T')} \left(
|\nabla(\Delta u(T^{\prime})|^2 + |\Delta u(T^{\prime})|^2 + 
\nabla u(T^{\prime}|^2  \right)\ dx.
\end{eqnarray*}

So we get for $s$ sufficiently large
\begin{eqnarray*}
s^2 \lambda^2 \int_{\Omega}  e^{-2s \eta(T')} \varphi(T')( |\nabla \gamma|^2\ + | \gamma|^2 ) \ dx 
\leq C 
    s \lambda  \int_{t_0}^T \int_{\Gamma_0} e^{-2s \eta} \varphi | \partial_{\nu} y |^2\ d x\ d t\\ 
   + C  \int_{\Omega}  e^{-2s \eta(T')} \left(
|\nabla(\Delta u(T^{\prime})|^2 + |\Delta u(T^{\prime})|^2 + 
\nabla u(T^{\prime}|^2  \right)\ dx,
\end{eqnarray*}

 and the the theorem is proved.
\end{proof}

{\bf Remark}\\
\begin{itemize}
\item
All the previous results are available for $\Omega \subset \mathbb{R}^n$ be a bounded domain of $\mathbb{R}^n$
with $n \geq 3$ if we adapt the regularity properties of the initial and boundary data.
\item
We give a stability result for two linked coefficient ($c$ and $\nabla c$) with one observation. 
Note that for two independent coefficients, there is no result in the litterature with only one
observation.
\end{itemize}

\medskip
 {\it E-mail address: }gaitan@cmi.univ-mrs.fr\\
 
\end{document}